\documentclass[12pt]{article}

\usepackage{amsfonts}     
\usepackage{amsmath}
\usepackage{amstext}
\usepackage{amsthm}
\usepackage{xspace}

\input psfig.sty

\newcommand{\N}{\mathbb N}

\newtheorem{teorema} {Theorem}[section]
\newtheorem{prop}[teorema]{Proposition}

\theoremstyle{definition}
\newtheorem{definizione}[teorema]{Definition}
\newtheorem{cor}[teorema]{Corollary}
\newtheorem{esempio}[teorema]{Example}
\newtheorem{guess}[teorema]{Remark}

\begin{document}                                            
\title{The Manneville map: topological, metric and algorithmic entropy}
\author{Claudio Bonanno \\
\small Dipartimento di Matematica \\ 
\small Universit\`a di Pisa \\ 
\small via Buonarroti 2, 56100 Pisa (Italy) \\ 
\small email: bonanno@mail.dm.unipi.it}
\date{}
\maketitle
\begin{abstract}
We study the Manneville map $f(x)=x+x^z (mod\ 1)$, with $z>1$, from a
computational point of view, studying the behaviour of the {\it
Algorithmic Information Content}. In particular, we consider a family
of piecewise linear maps that gives examples of algorithmic behaviour
ranging from the {\it fully} to the {\it mildly} chaotic, and show
that the Manneville map is a member of this family.
\end{abstract}
\section{Introduction} \label{sI}
\begin{figure}[h]
\psfig{figure=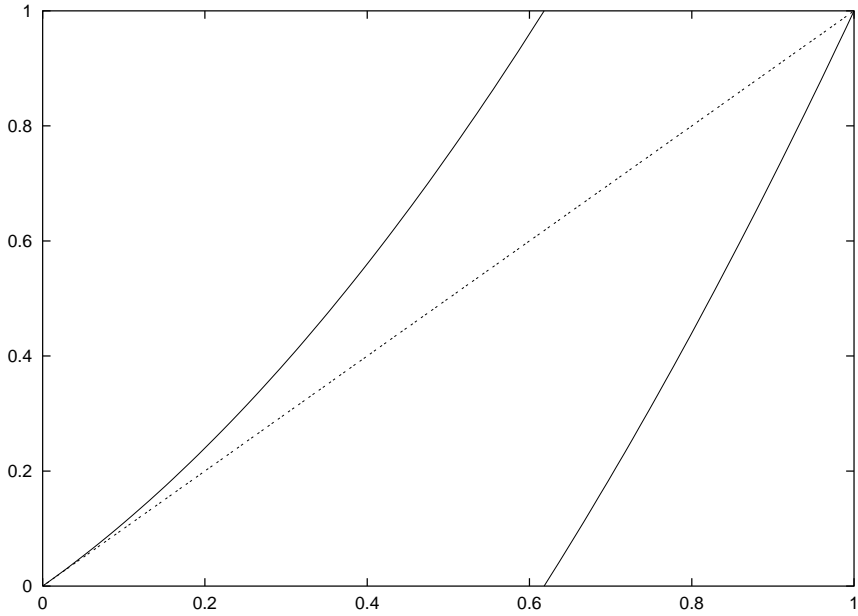,width=12cm,angle=270}
\caption{\it The Manneville map $f$ for $z=2$}
\label{figmann}
\end{figure}

The {\it Manneville map} was introduced by Manneville in \cite{Mann},
as an example of a discrete dissipative dynamical system with {\it
intermittency}, an alternation between long regular phases, called
{\it laminar}, and short irregular phases, called {\it
turbulent}. This behaviour had been observed in fluid dynamics
experiments and in chemical reactions. Manneville introduced his map,
defined on the interval $I=[0,1]$ by
\begin{equation}
f(x)=x+x^z (\hbox{mod } 1) \hskip 0.5cm z > 1, \label{mannmap}
\end{equation}
to have a simple model displaying this complicated behaviour. In
Figure \ref{figmann} it is plotted the Manneville map for $z=2$. His
work has attracted much attention, and the dynamics of the Manneville
map has been found in many other systems. We can find applications of
the Manneville map in dynamical approaches to DNA sequences
(\cite{Grigo1},\cite{Grigo2}) and ion channels (\cite{Toth}), and in
non-extensive thermodynamical problems (\cite{Grigo3}).

The Manneville map has also been studied by Gaspard and Wang
(\cite{GW}), using the notion of {\it Algorithmic Information Content}
of a string, briefly explained below. Given a dynamical system, any
orbit of this system can be translated into a string $\sigma$ of
symbols by a partition of the phase space of the system (symbolic
dynamics). For any finite string $\sigma^n$ of length $n$, it has been
introduced by Chaitin (\cite{Chaitin}) and Kolmogorov (\cite{Kolm}),
the notion of {\it Algorithmic Information Content (AIC)} (or {\it
Kolmogorov complexity}) of the string, that we denote by
$I_{AIC}(\sigma^n)$, defined as the binary length of the shortest
program $p$ that outputs the string on a universal machine $C$,
\begin{equation}
I_{AIC}(\sigma^n)= \left\{ |p| \ | \ C(p)= \sigma^n \right\}. 
\label{AIC} 
\end{equation}

It is then possible, using the symbolic dynamics, to define the notion
of Algorithmic Information Content for a finite orbit of a dynamical
system. This extension requires some attention on the choice of the
partition. The first results have been obtained by Brudno
(\cite{Brudno}) using open covers of the phase space. Another possible
approach, using computable partitions, is introduced in \cite{Licata}.

To generalize the notion of AIC to infinite strings, it is natural to
consider the mean of the AIC. We call {\it complexity} of an infinite
string $\sigma$ the maximum limit of the AIC of the first $n$ symbols
of the string divided by $n$. Then, if we denote the complexity of an
infinite string by $K(\sigma)$, we have
\begin{equation}
K(\sigma)=\limsup_{n\to +\infty} \frac{I_{AIC}(\sigma^n)}{n},
\label{complessita}
\end{equation}
where $\sigma^n$ is the string given by the first $n$ digits of the
infinite string $\sigma$. Symbolic dynamics is again the tool to
define the complexity of an infinite orbit of a dynamical system.

Moreover, we can ask whether it is possible to define a notion of
information content for the dynamical system, without to consider any
particular orbit. To do this, we have to introduce a probability
measure $\mu$ on the phase space $X$ of the system and we can define
the {\it algorithmic entropy} $h_\mu$ of a dynamical system by
\begin{equation}
h_\mu = \int_X \ K(x) \ d\mu, \label{entropia}
\end{equation}
where $K(x)$ denotes the complexity of the orbit of the system with
initial condition $x \in X$.

There exist some results connecting the information content of a
string generated by a dynamical system and the Kolmogorov-Sinai
entropy $h^{KS}$ of the system. 

First of all it is proved that for a compact phase space $X$ and for
an invariant measure $\mu$, we have $h_\mu=h^{KS}_\mu$. Then in
particular, in a dynamical system with an ergodic invariant measure
$\mu$ with positive K-S entropy $h^{KS}_\mu$, the AIC of a string $n$
symbols long behaves like $I_{AIC}(\sigma^n) \sim h^{KS}_\mu n$ for
almost any initial condition with respect to the measure $\mu$
(\cite{Brudno}).

Instead, in a periodic dynamical system, we expect to find
$I_{AIC}(\sigma^n) = O(\log(n))$. Indeed, the shortest program that
outputs the string $\sigma^n$ would contain only information on the
period of the string and on its length.

It is possible to have also intermediate cases, in which the K-S
entropy is null for all the invariant measures that are physically
relevant and the system is not periodic. These systems, whose
behaviour has been defined {\it weak chaos}, are an important
challenge for research on dynamical systems. Indeed no information are
given by the classical properties, such as K-S entropy or Lyapunov
exponents, and in the last years some generalized definitions of
entropy of a system have been introduced to characterize the behaviour
of such systems (for example see \cite{Tsallis}). We believe that an
approach to weakly chaotic systems using the infinite order of their
AIC could be a powerful way to classify these systems (no information
are obtained by the complexity and the algorithmic entropy defined as
above).

The Manneville map with parameter $z>2$ is a non periodic map with
null K-S entropy for all the physically relevant invariant measures,
then the analysis of the AIC of the strings generated by the map is
interesting. Gaspard and Wang (\cite{GW}) showed that the Manneville
map exhibits a behaviour that they called {\it sporadicity}. Namely,
the Algorithmic Information Content, $I_{AIC}(\sigma^n)$, of a string
$n$ symbols long, behaves in mean like $n^\alpha (\log(n))^\beta $,
with either $0<\alpha <1$ or $\alpha =1$ and $\beta <0$.

In this paper, we give a formal proof of the results obtained by
Gaspard and Wang (\cite{GW}) for the Manneville map, giving more
precise estimates for the AIC of a string generated by the map. But
the most important generalization is that we find our results for the
Manneville map as a particular case of a general theorem concerning a
large family of maps $L$, defined in equation (\ref{mannmaplin}). This
family of maps exhibits an extremely wide range of behaviours (for the
AIC of the generated strings), and sporadicity is only one possible
case. Then we find a family of maps that can be classified with
respect to the order of the AIC of a ``typical''(in the sense of
Lebesgue measure) initial condition.  Moreover we study some
topological and metric properties of the maps $L$, useful to obtain a
prediction of the behaviour of the AIC of the related symbolic
dynamics (see Section \ref{sTca}).

In Section \ref{sTlm}, we introduce the family of maps $L$, and show
how the Manneville map $f$ can be thought of as one member of the
family.

In Section \ref{sTta}, we study the maps $L$ from the topological
point of view. In Section \ref{sTma} we show that the maps in the
class $L$ are equivalent to a Markov chain in a suitable sense. This
equivalence is extensively used in Subsection \ref{ssGr}, where we
present our results relative to the behaviour of the AIC of the
strings obtained from the maps $L$.

Finally, in Subsection \ref{ssRttlm}, we study the computational
aspect of the maps $L$ from a practical point of view. So, we restrict
ourselves to consider the Lebesgue measure $l$ on the interval $I$.

\section{The family of piecewise linear maps $L$} \label{sTlm}
In this section we present what we shall use as our formulation of the
Manneville map (\ref{mannmap}). In the following, we study a family of
piecewise linear maps $L$ on the interval $I=[0,1]$, which are
topologically equivalent to the Manneville map $f$. Using the maps
$L$, all the theorems have an easier interpretation and computations
can be done exactly. Moreover all our results are extendible through
metric isomorphism, hence we shall define on the interval $I$, two
different measures that make the topological equivalence between $L$
and $f$ a metric isomorphism. Then we can extend all the results that
we find for the maps $L$ to the Manneville map $f$.

Let's start defining the piecewise linear maps $L$ that we consider. We
use here the same approach as in \cite{Stefano}. A natural way to get a
partition of the interval $I=[0,1]$ from the Manneville map $f$ is the
following: let's call $x_0$ the point of $I$ such that $f(x_0)=1$ with
$x_0 \not= 0,1$, and $x_1$ the preimage of $x_0$ in the interval
$[0,x_0]$; then we define recursively $x_n = \{ f^{-1}(x_{n-1}) \}
\cap [0,x_{n-1}]$. Then the sub-intervals $B_k = (x_k, x_{k-1}]$, for
$k \geq 1$, and $B_0 = (x_0,1]$ are a partition of $I$.

Define $\{ \epsilon_k \}_{k\in \N}$ a sequence of positive real
numbers, that is strictly monotonically decreasing and converging
towards zero, with the property that
\begin{equation}
\frac{\epsilon_{k-1}-\epsilon_k}{\epsilon_{k-2}-\epsilon_{k-1}} < 1
\hskip 0.5cm \forall \ k \in \N. \label{propdieps}
\end{equation}

The piecewise linear maps that we consider are defined by
\begin{equation}
L(x)=\left\{
\begin{array}{ll}
\frac{\epsilon_{k-2}-\epsilon_{k-1}}{\epsilon_{k-1}-\epsilon_k} (x -
\epsilon_k) + \epsilon_{k-1} & \quad \epsilon_k < x \leq
\epsilon_{k-1}, \quad k \geq 1 \\[0.3cm]
\frac{x-\epsilon_0}{1-\epsilon_0} & \quad \epsilon_0 < x \leq 1
\\[0.3cm] 
0 & \quad x=0
\end{array}
\right. \label{mannmaplin}
\end{equation}
where we define $\epsilon_{-1} = 1$. These piecewise linear maps $L$
clearly depend on the definition of the sequence $\{ \epsilon_k
\}_{k\in \N}$, but a particular choice for this sequence is important
only for Section \ref{sTca}. For the moment we consider any possible
sequence, with the properties specified above. Let's define
sub-intervals $A_i = (\epsilon_i, \epsilon_{i-1}]$, and
$A_0=(\epsilon_0,1]$. These interval form a partition of the interval
$I$. We prove the following

\begin{teorema}
Any piecewise linear map $L$ defined as in equation (\ref{mannmaplin})
is topologically equivalent to the Manneville map $f$. \label{teoeqlin}
\end{teorema}

\noindent {\bf Proof.} We have to find a homeomorphism $h :I \to I$
such that $h(f(x))=L(h(x))$ for each $x \in I$. To find such a
homeomorphism we use the partitions $(B_j)$ and $(A_j)$. We define
$h(x_n)=\epsilon_n$ for each $n \in \N$ and $h(0)=0$, $h(1)=1$, and
such that $h(B_j)=A_j$, for all $j \geq 0$. To define the
homeomorphism $h$ we use a dense set of $I$, define $h$ on this set
and, then, simply extends the definition of $h$ to the whole interval
$I$ by continuity.

Let's consider a sub-interval $B_k$. By definition of the Manneville
map $f$, we have that $f(B_k)=B_{k-1}$, for $k\geq 1$, and
$f(B_0)=I$. Then it follows that $f^k(B_k)=B_0$ and
$f^{k+1}(B_k)=I$. So, within each $B_k$ we can find sub-intervals
$B_{kj}$, with $k,j \in \N$, defined by $f^{k+1}(B_{kj})=B_j$. These
sub-intervals form a partition of each $B_k$. We can continue this
partition of the intervals $B_k$, defining, by the same rule,
sub-intervals $B_{kji}$ that form a partition of $B_{kj}$. We write
then any sub-interval of the form $B_{k_1 k_2 \dots k_n}$, with $k_i
\in \N$ for each $i=1,\dots ,n$, as $B_{k_1 k_2 \dots k_n}=(x_{k_1 k_2
\dots k_n}, x_{k_1 k_2 \dots (k_n-1)}]$. The set $\left\{ x_{k_1 k_2
\dots k_n}, \ n\in \N, \ k_i \in \N \right\}$ is a countable dense set
of the interval $I$. Analogously, we can define a set of points
$\left\{ \epsilon_{k_1 k_2 \dots k_n}, \ n\in \N, \ k_i \in \N
\right\}$, for the map $L$, with the same property as $x_{k_1
k_2 \dots k_n}$. We define then $h(x_{k_1 k_2 \dots
k_n})=\epsilon_{k_1 k_2 \dots k_n}$ for each $n\in \N$, and extend the
function $h$ continuously to the whole interval $I$. We have thus
obtained a continuous function $h$ such that $h(B_{k_1 k_2 \dots
k_n})=A_{k_1 k_2 \dots k_n}$, for each $n\in \N$, where the
sub-intervals $A_{k_1 k_2 \dots k_n}$ are defined as $A_{k_1 k_2 \dots
k_n}=(\epsilon_{k_1 k_2 \dots k_n}, \epsilon_{k_1 k_2 \dots
(k_n-1)}]$.

The injectivity of $h$ follows by contradiction. Let's suppose to have
two points $x < y \in I$ such that $h(x)=h(y)=z$. By the density of
the set $\left\{ x_{k_1 k_2 \dots k_n}, \ n\in \N, \ k_i \in \N
\right\}$, we can find a point $\bar x_{k_1 k_2 \dots k_{\bar n}} \in
(x,y)$. This implies that there exists a $\bar n$ such that $x \in
B_{k_1 k_2 \dots (k_{\bar n}-1)}$ and $y \in B_{k_1 k_2 \dots k_{\bar
n}}$. Then we have $h(x)\not= h(y)$. The subjectivity of $h$ follows
immediately by the definition.

Then the inverse function $h^{-1}$ exists and $h$ is a homeomorphism
because it is a continuous invertible function from a compact to a
Hausdorff space. \qed

\vskip 0.5cm The topological equivalence between $L$ and $f$ can be
used in particular to obtain a metric isomorphism. If we have a
measure $\mu$ on the interval $(I,{\cal B},L)$, where ${\cal B}$ is
the Borel $\sigma$-algebra, then the homeomorphism $h$ carries $\mu$
into another measure $\nu=h^* \mu$ on $(I,{\cal B},f)$, and with
respect to these measures $h$ is a metric isomorphism. 

\begin{teorema}[Radon-Nikodym]
Given two measures $\mu$ and $\nu$ on $(I,{\cal B})$, such that
$(I,{\cal B},\nu)$ is $\sigma$-finite, $\mu << \nu$ if and only if
there exists a real function $f$ on $I$, integrable with respect to
$\nu$ on all sets $B \in {\cal B}$ such that $\nu (B) < +\infty$,
satisfying the following condition for every $B \in {\cal B}$:
\[
\mu(B) = \int_B f \ d\nu
\]
\label{rdteo}
\end{teorema}

\begin{teorema}
If the measure $\mu$ on $(I,{\cal B},L)$ is absolutely continuous with
respect to the Lebesgue measure $l$ ($\mu << l$), then also the
measure $\nu = h^* \mu$ on $(I,{\cal B},f)$ is absolutely continuous
with respect to the Lebesgue measure $l$, and vice-versa. \label{miseq}
\end{teorema}

\noindent {\bf Proof.} We can apply Theorem \ref{rdteo} to the
measures $\mu$ and $l$. Then we obtain a real function $f_\mu$ on $I$,
integrable with respect to $l$ and such that for all $B \in {\cal B}$
\[
\mu(B)= \int_B f_\mu \ dl.
\]

By definition of the measure $\nu$, we have that $\nu(B) = \mu(h(B))$
for all $B \in {\cal B}$, then
\[
\nu(B)= \int_{h(B)} f_\mu \ dl = \int_B (f_\mu \circ h) (dh) \ dl,
\]
where $dh$ is defined almost everywhere with respect to $l$, being $h$
a monotone continuous function. Then we have found a function $f_\nu =
(f_\mu \circ h) (dh)$, which satisfies the hypotheses of the
Radon-Nikodym Theorem. Then $\nu << l$.

The vice-versa is proved in the same way. \qed

\vskip 0.5cm At this point, thanks to Theorem \ref{miseq}, we can use
our linear map in all our applications of topological and metric
methods. 

\section{The topological approach} \label{sTta}
In this section we start a procedure of equivalences of the maps $L$
with well-known maps, that can be used to establish the results of
Section \ref{sTca}. The first step is to study the relationship
between the maps $L$ and a {\it sub-shift of finite type}.

\subsection{Symbolic dynamics} \label{ssSd}
We briefly recall the definition of a {\it sub-shift of finite
type}. Let's consider a finite set of symbols $S= \{ 0,1,2,\dots,N\}$,
with $N \geq 1$, and build a set $\Sigma^N = S^{\N}$ as the product of
countable factors $S$. On $\Sigma^N$ it is defined a map $T$, called
the {\it shift map}, that acts on the elements of $\Sigma^N$, by
shifting forward the indexes. Namely, if $\sigma = (\sigma_0 \sigma_1
\dots \sigma_n \dots) \in \Sigma^N$, with $\sigma_i \in \{ 0,\dots,N
\}$ for all $i \in \N$, then $T(\sigma)= (\sigma_1 \dots \sigma_n
\dots)$. The set $\Sigma^N$ is endowed with a metric $d$ defined by
\[
d(\sigma,\sigma')= \sum_{n=0}^{\infty} \ \frac{\delta_{\sigma_n
\sigma'_n}}{2^n},
\]
where $\delta_{ij}$ is the Kronecker symbol, that makes it a compact
space. A {\it sub-shift of finite type} is obtained from the set
$(\Sigma^N,T)$, by means of a $(N+1)\times (N+1)$ matrix $M=(m_{ij})$, called
the {\it transition matrix}, such that $m_{ij} \in \{0,1 \}$ for all
$i,j=0,\dots, N$. We define a subset $\Sigma^N_M$ of $\Sigma^N$, by
\[
\Sigma^N_M = \left\{ \sigma \in \Sigma^N \ | \ m_{\sigma_i
\sigma_{i+1}} =1 \; \forall \; i \in \N \right\},
\]
then a {\it sub-shift of finite type} is simply the compact, metric
space $(\Sigma^N_M,T_M)$ (see \cite{KH}).

We have 

\begin{teorema}
For any $N\geq 1$, there exists a particular transition matrix $M$
such that any map $L$ of the family of piecewise linear maps
(\ref{mannmaplin}) is a factor of the sub-shift of finite type
$(\Sigma^N_M,T_M)$. This means that there exists a subjective
continuous map $\pi: \Sigma^N_M \to I$ such that $\pi \circ T_M = L
\circ \pi$. \label{eqdinsimb}
\end{teorema}

\noindent {\bf Proof.} Let a $(N+1)\times (N+1)$ matrix $M$ be defined
by
\begin{equation}
M= \left(
\begin{array}{ccccccc}
1 & 1 & 1 & \cdots & 1 & 1 & 1\\
1 & 0 & 0 & \cdots & 0 & 0 & 0\\
0 & 1 & 0 & \cdots & 0 & 0 & 0\\
\cdots & \cdots & \cdots & \cdots & \cdots & \cdots & \cdots\\
0 & 0 & 0 & \cdots & 0 & 1 & 1
\end{array}
\right), \label{transmatrix}
\end{equation}

\noindent we show that this is the matrix we need.

Let's consider a partition of the interval $I$ defined by using the
partition $A_i =(\epsilon_i , \epsilon_{i-1}]$, introduced in Section
\ref{sTlm}. From this partition we obtain a partition of the interval
$I$, by $B_i = A_i$ for $i=0,\dots, N-1$, and $B_N = \{ 0 \} \cup
\left( \cup_{j=N}^{\infty} A_j \right)= [0,\epsilon_{N-1})$. From this
partition and by using the properties of any of the map $L$, we obtain
a $(N+1)$-nary representation of the interval $I$.

The $(N+1)$-nary representation of a point $x\in I$ is given by a
string $\sigma$ such that $L^n(x) \in B_{\sigma_n}$ with $\sigma_n
=0,\dots,N$, for any $n\in \N$. This representation is nothing else
that a map $\pi : \Sigma^N \to I$. Hence we just need to show that
this map $\pi$ is continuous and subjective, and verifies the
commutation rule with $L$ and $T$.

First of all we notice that our map $\pi$ is not defined on the whole
space $\Sigma^N$, because of the restrictions given by the particular
form of the map $L$. If we want to reduce the space $\Sigma^N$, we
have to consider a transition matrix $M$, and we use the matrix $M$
defined in equation (\ref{transmatrix}). It is easy to verify that for
any $\sigma \in \Sigma^N_M$ there is a point $x\in I$ such that
$\pi(\sigma)=x$. Then we show that $\pi : \Sigma^N_M \to I$ is a
subjective continuous map such that $\pi \circ T_M = L \circ \pi$. 

We have then a {\it semi-conjugacy} between our piecewise linear maps
$L$ and symbolic dynamics. We have not a conjugacy because of the lack
of injectivity of the map $\pi$. Indeed, as in any $n$-ary
representation of the real numbers in the interval $I=[0,1]$, there is
a countable set $X$ of points that are images of two
sequences. Moreover, in our case, these points can be characterized by
the property that for any $x$ in this set there exists a $N \in \N$
such that $L^N(x)=1$. 

Commutation. The commutation rule $\pi \circ T_M = L \circ \pi$ is an
immediate consequence of the definition of $\pi$.

Subjectivity. It follows immediately from the definition of the map $\pi$.

Continuity. We have to prove that given any $\epsilon >0$ there exists
a $\delta >0$ such that if $d(\sigma_1,\sigma_2) <\delta$ then
$|\pi(\sigma^1)-\pi(\sigma^2)| <\epsilon$. But from the definition of
the metric $d$ on the space $\Sigma^N_M$, we have that
$d(\sigma^1,\sigma^2) <\delta$ is equivalent to: there exists a $K>0$
such that $\sigma^1_j=\sigma^2_j$ for all $j=0,\dots,K$. So given any
$\epsilon>0$ we have to find a $K>0$ such that $\sigma^1_j=\sigma^2_j$
for all $j=0,\dots,K$ implies $|\pi(\sigma^1)-\pi(\sigma^2)|
<\epsilon$. From the definition of the map $\pi$ it is clear that if
$\sigma^1_j=\sigma^2_j$ for all $j=0,\dots,K$ for any $K>0$, then
$L^j(\pi(\sigma^1))$ and $L^j(\pi(\sigma_2))$ belong to the same
subset $B_{\sigma^i_j}$ of the partition $(B_i)$ for all
$j=0,\dots,K$. Then, if we consider a partition of the subset
$B_{\sigma^i_0}$, given by $(B_{\sigma^i_0})_{j_1,j_2,\dots,j_n}$,
with $j_i=0,\dots,N$ for all $i$, where
$L^r((B_{\sigma^i_0})_{j_1,j_2,\dots,j_n}) = B_{j_r}$ for all
$r=1,\dots,n$, we have that
$diam((B_{\sigma^i_0})_{j_1,j_2,\dots,j_n}) \to 0$ as $n \to +\infty$,
thanks to the particular form of the map $L$. This argument gives the
continuity of the map $\pi$.\qed

\vskip 0.5cm

We can extend Theorem \ref{eqdinsimb} to the case of $N=\infty$. The
space $\Sigma=\Sigma^{\infty}$ is defined in the same way, but we
cannot extend the metric $d$, defined as before, and we can only
define a topology on $\Sigma$, where the open balls are the same as
before. In this case $\Sigma$ is not anymore a compact space. The
transition matrix $M$ is $\infty \times \infty$ dimensional and it is
defined by
\begin{equation}
M= \left(
\begin{array}{cccccccc}
1 & 1 & 1 & \cdots & 1 & 1 & 1 & \cdots \\
1 & 0 & 0 & \cdots & 0 & 0 & 0 & \cdots \\
0 & 1 & 0 & \cdots & 0 & 0 & 0 & \cdots \\
0 & 0 & 1 & \cdots & 0 & 0 & 0 & \cdots \\
\cdots & \cdots & \cdots & \cdots & \cdots & \cdots & \cdots & \cdots
\end{array}
\right). \label{inftransmatrix}
\end{equation}

On the space $\Sigma_M$, we define then a map $T_M$, given by the
forward shift. Then we can prove 

\begin{teorema}
Any of our maps $L$ on $I$ is a topologically conjugate to the
dynamical system given by $(\Sigma_M,T_M)$ built on countable symbols
and transition matrix $M$ given by equation (\ref{inftransmatrix}).
\label{eqinfdinsimb}
\end{teorema}

\noindent {\bf Proof.} The proof is the same as in Theorem
\ref{eqdinsimb}. The difference is that now we obtain a {\it
conjugacy} with the set $\Sigma_M$, thanks to the fact that the matrix
$M$, being infinite, exactly simulate the dynamics of the piecewise
linear maps $L$ from a topological point of view. \qed

\subsection{Topological and Kolmogorov-Sinai entropy} \label{ssTame}
The semi-conjugacy of the maps $L$ with the symbolic dynamical system
on finite symbols is useful to compute some topological and metric
quantities of our map. Indeed the lack of injectivity of the map $\pi$
is on a set that doesn't change the dynamical richness of the
systems. In this subsection we compute the {\it topological entropy}
and the {\it Kolmogorov-Sinai entropy} for some measures (see
\cite{KH}).

From the theory of dynamical systems, we know that the following
theorems hold (see \cite{KH}):

\begin{teorema}
The topological entropy is invariant for topological
equivalence. \label{invarianzaet}
\end{teorema}

\begin{teorema}
The topological entropy $h_{top}$ of a sub-shift of finite type is
$\log \lambda_{\max}$, where $\lambda_{\max}$ is the largest
eigenvalue of the transition matrix. \label{tedinsimb}
\end{teorema}

Then we have just to compute the eigenvalues of the transition matrix
$M$ on finite symbols defined in equation (\ref{transmatrix}), and
then, thanks to the previous theorems, the topological entropy
$h_{top} (L)$ of our linear map $L$ is given by $\log
\lambda_{\max}$. For any $N$, we find that $\lambda_{\max} =2$, then
$h_{top} (L)= \log 2$, for any possible sequence $(\epsilon_k)$
defined as before.

At this point we start to consider measures on the space
$\Sigma^N_M$. We have to introduce first a $\sigma$-algebra ${\cal
C}$. We take as a basis of $\cal C$ the sets of the form
\begin{equation}
C^n_r = \left\{ \sigma \in \Sigma^N_M \ | \ \sigma_i = r_i \; \forall \; i =0,\dots n \right\}, \label{cilindri}
\end{equation}
for any $n \in \N$ and $r \in S^{n+1}$. These sets are called {\it
cylinders}. At this point we use a classical result of dynamical
systems (see \cite{Mane}):

\begin{teorema}[Variational Principle]
Given a continuous map $f: X \to X$ of a compact metric space $X$, the
topological entropy $h_{top}(f)$ is the maximum of the
K-S entropies $h_{\nu} (f)$ on the set of all the
$f$-invariant probability measures $\nu$ on $X$.
\end{teorema}

Then we look for the probability measures $\nu$ on $\Sigma^N_M$ with
K-S entropy $h^{KS}_\nu$ equal to $\log 2$. For sub-shift of finite
type, a particular class of $T_M$-invariant measures are defined by a
{\it stochastic matrix} associated to the sub-shift. These measures
are called {\it Markov measures}, and among them there is a measure
that maximize the K-S entropy.  This measure is called {\it Parry
measure}, and is denoted by $\nu_\Pi$ (see \cite{KH}). This measure is
defined by a particular choice of the stochastic matrix $\Pi$. In
words, the Parry measure represent the asymptotic distribution of the
periodic orbits, that is if $C$ is a cylinder in $\Sigma^N_M$, we have
that
\[
\nu_\Pi (C) = \lim_{n \to \infty} \frac{\hbox{periodic orbits of
period $n$ contained in $C$}}{\hbox{all the periodic orbits of period
$n$}}.
\]

We compute the Parry measures for some $N$, then we define on $I$ the
induced measure $\mu_\Pi=\pi^* \nu_\Pi$, to obtain an $L$-invariant
probability measure on $I$ of K-S entropy
$h^{KS}_{\mu_\Pi}(L)= \log 2$.

For $N=1$, we have that $\nu_\Pi (C^0_0) = \nu_\Pi (C^0_1) =
\frac{1}{2}$. So $\mu_\Pi (B_0) = \mu_\Pi (B_1) = \frac{1}{2}$.
For any $N$, we obtain
\[
\begin{array}{ll}
\mu_\Pi (B_i) & = \frac{1}{2^{i+1}} \qquad i=0,\dots,N-1 \\[0.3cm]
\mu_\Pi (B_N) & = \frac{1}{2^{N-1}}
\end{array}
\]

So for the limit $N \to \infty$ (where we directly have topological
equivalence between $(I,L)$ and $(\Sigma_M,T_M)$), we obtain $\mu_\Pi
(B_i) = \mu_\Pi (A_i) = \frac{1}{2^{i+1}}$. We have thus found a
countable family of $L$-invariant measures on $I$ with K-S entropy
$\log 2$.

\section{The metric approach} \label{sTma}
We start now to consider what happens if we start from our interval
$I$, endowed with a probability measure $\mu$ on the Borel
$\sigma$-algebra $\cal B$, and with the dynamics induced by the maps
$L$. In particular we show that we obtain an equivalence with a Markov
chain that will be useful for the computational approach (Section
\ref{sTca}).

From our space $(I,{\cal B},L,\mu)$, where we remark that we haven't
supposed $\mu$ to be $L$-invariant, we have a metric isomorphism with
the space $(\Sigma_M,{\cal C},T_m)$ on countable symbols, with the
probability measure $P=\pi_* \mu$ induced by the homeomorphism $\pi$
found in Theorem \ref{eqinfdinsimb}. 

At this point we use some notions and results introduced by Parry
\cite{Parry}. 

\begin{definizione}
A {\it non-atomic stochastic process} is $(X,{\cal A},T,m)$ where
$X=\{ x= x_0, x_1, \dots \; | \; x_i \in \N \; \forall i \in \N \}$,
${\cal A}$ is the $\sigma$-algebra generated by the cylinders $C^n_r$,
$m$ is a non-atomic probability measure on ${\cal A}$, and $T$ is the
forward shift on $X$. In the theory of stochastic processes the
transition matrix $M$ on $\Sigma_M$ is called a {\it structure matrix}.
\end{definizione}

\begin{definizione}
A stochastic process is called {\it transitive of order $k$} if for
all $(x_1, \dots, x_k)$ and $(y_1,\dots,y_k)$ with $m(x_1, \dots, x_k)
>0$ and $m(y_1, \dots, y_k) >0$, there exists a finite
$(z_1,\dots,z_n)$ such that 
\[
m(x_1, \dots, x_k;z_1, \dots, z_n;y_1, \dots, y_k) >0.
\]
When $k=1$, a stochastic process is simply called {\it transitive}.
\end{definizione}

\begin{definizione}
A stochastic process is said to be {\it intrinsically Markovian of
order $k$} if $m(x_1,\dots,x_n)>0$ and $m(x_{n-k+1},\dots,x_{n+1})>0$
imply 
\[
m(x_1,\dots,x_n,x_{n+1})>0. 
\]
When $k=1$ it is simply called {\it intrinsically Markovian}.
\end{definizione}

We can easily prove

\begin{prop}
Our space $(\Sigma_M,{\cal C},T_M,P)$ is a non-atomic stochastic
process, which is transitive and intrinsically Markovian. \label{prop1}
\end{prop}

\begin{definizione}
Given a stochastic process $(X,{\cal A},T,m)$, a measure $p$ makes the
process $(X,{\cal A},T,p)$ {\it compatible} with the original when
$p(C^n_r)>0$ if and only if $m(C^n_r)>0$, for any cylinder $C^n_r$.
\end{definizione}

At this point we use the notion of non-atomic stochastic processes to
obtain a compatibility between the maps $L$ and a Markov chain,
through the symbolic dynamical system on countable symbols. Before
giving the theorems in this direction, we briefly recall the theory of
Markov chains (see \cite{Chung}).

\begin{definizione}
Given a probability space $(\Lambda, {\cal F}, P)$ and a countable
space $Y$ with the discrete $\sigma$-algebra, a {\it Markov chain} is
a sequence $(Z_n)_{n\in \N}$ of random variables $Z_n : \Lambda \to Y$ such that

i) If, given $y_0,\dots,y_{n+1} \in Y$, we have $P[Z_n=y_n,
Z_{n-1}=y_{n-1},\dots,Z_0=y_0]>0$, then
\[
P[Z_{n+1}=y_{n+1}\; | \; Z_n=y_n,\dots,Z_0=y_0] = P[Z_{n+1}=y_{n+1}\; | \; Z_n=y_n],
\]

ii) If $x,y \in Y$ and $m,n \in \N$ are such that $P[Z_m=x]>0$ and
$P[Z_n=x]>0$, then
\[
P[Z_{m+1}=y \; | \; Z_m=x]=P[Z_{n+1}=y \; | \; Z_n=x].
\]

In particular the numbers $p(x,y)=P[Z_{n+1}=y \; | \; Z_n=x]$ form a
matrix $\Pi = (p(x,y))_{x,y \in X}$, called the {\it transition
matrix}. Moreover the probability measure on $Y$ defined by $\nu(y) =
P[Z_0=y]$ is called the {\it initial distribution}.
\end{definizione}

The transition matrix $\Pi$ is a {\it stochastic matrix}, that is
$p(x,y)\geq 0$ and $\sum_{y \in X} p(x,y) =1$ for all $x \in Y$.

\begin{teorema}
Given any countable space $Y$, a transition matrix $\Pi$ and an
initial distribution $\nu$, it is possible to construct a probability
space $(\Lambda,{\cal F},P)$ and a sequence of random variables
$Z_n:\Lambda \to Y$, such that the constructed Markov chain has $\Pi$
as transition matrix and $\nu$ as initial distribution. \label{yonescotulcea}
\end{teorema}

\noindent {\bf Proof.} For the proof of the theorem see Chung
\cite{Chung}. We simply say what is the constructed probability space
$(\Lambda,{\cal F},P)$. The space $\Lambda$ is $Y^{\N}$ and is called
the {\it realizations space}, the $\sigma$-algebra ${\cal F}$ is given
by the cylinders defined as in equation (\ref{cilindri}), and the
probability $P$ is defined on the cylinders by
\[
P(C^n_r)=\nu(r_0) p(r_0, r_1) p(r_1, r_2) \dots p(r_{n-1} r_n).
\]
The random variables are defined as the projections of $\Lambda$ on $Y$.
\qed

\begin{teorema}
Given any intrinsically Markovian, transitive sto\-cha\-stic process
$(X,{\cal A},T,m)$ and a stochastic matrix $\Pi$ such that $p(i,j)>0$
if and only if $m_{ij}=1$, for the structure matrix $M$ of the
process, there is a probability $p$ on $X$ such that $(X,{\cal
A},T,p)$ is compatible with $(X,{\cal A},T,m)$, and it is a Markov
chain with $\Pi$ as transition matrix. \label{compatibili}
\end{teorema}

\noindent {\bf Proof.}
We have just to apply Theorem \ref{yonescotulcea} to the matrix $\Pi$
and to an initial distribution $\nu$, being $X$ already in the form of
the realizations space. The probability $p$ is then the probability
$P$ defined as above.
\qed

\begin{cor}
Our space $(\Sigma_M,{\cal C},T_M,P)$ is com\-pa\-ti\-ble with a Mar\-kov
cha\-in. \label{cor1}
\end{cor}

\noindent {\bf Proof.}  The corollary is proved thanks to Theorem
\ref{compatibili} and Proposition \ref{prop1}. We have just to choose
a stochastic matrix that satisfies the hypothesis, that is $p(i,j) >0
$ if and only if $m_{ij}=1$ for the structure matrix defined as in
equation (\ref{inftransmatrix}). The measure $P$ can be used as
initial distribution.
\qed

\vskip 0.5cm We have thus completed our equivalence, in the sense of
Corollary \ref{cor1}, between the maps $L$ defined on the space
$(I,{\cal B},\mu)$, for any probability measure $\mu$, and a Markov
chain, that is denoted simply by a stochastic matrix $\Pi$ and an
initial distribution $\nu$.

\section{The algorithmic entropy} \label{sTca}
The results of Section \ref{sTma} are useful for the computational
approach to the Manneville map $f$ defined by equation
(\ref{mannmap}). As remarked before, also in this section we shall
restrict ourselves to the AIC for the maps $L$, which are equivalent
to the Manneville map $f$ in the sense described above. Using this
restriction it is possible to perform explicit computations which, by
Theorems \ref{teoeqlin} and \ref{miseq}, can be extended to the
Manneville map $f$. The investigation on the maps $L$ that we present
in this section is meant to be a generalization of the work of Gaspard
and Wang on the Manneville map (see \cite{GW}).

Given our dynamical system $(I,{\cal B},L,\mu)$, where $\mu$ is any
probability measure on the Borel $\sigma$-algebra of the interval $I$,
we translate the orbit of a point $x \in I$ into a string
$\sigma=\pi(x) \in \Sigma_M$, with the transition matrix $M$ given in
equation (\ref{inftransmatrix}), and study the AIC of the
string. 

\subsection{General results} \label{ssGr}
In Section \ref{sTma}, we proved that for any map $L$ our dynamical
system $(I,{\cal B},L,\mu)$ is equivalent, in a sense specified above,
with a Markov chain with a stochastic matrix $\Pi$ defined by means of
the transition matrix $M$ of equation (\ref{inftransmatrix}), and a
given initial distribution that can be considered to be the measure
$\mu$ itself (see Corollary \ref{cor1}). We remark that it is not
necessary to choose the measure $\mu$ on $I$ to be $L$-invariant. Now,
we want to relate the dynamics of the Markov chain with our dynamical
system. A natural choice for the stochastic matrix $\Pi$ is
\begin{equation}
\Pi= \left(
\begin{array}{cccccccc}
p_0 & p_1 & p_2 & \cdots & p_n & p_{n+1} & p_{n+2} & \cdots \\
1 & 0 & 0 & \cdots & 0 & 0 & 0 & \cdots \\
0 & 1 & 0 & \cdots & 0 & 0 & 0 & \cdots \\
0 & 0 & 1 & \cdots & 0 & 0 & 0 & \cdots \\
\cdots & \cdots & \cdots & \cdots & \cdots & \cdots & \cdots & \cdots
\end{array}
\right), \label{stochmatrix}
\end{equation}
where $p_i$ are the probabilities of transition from the sub-interval
$A_0$ to the sub-interval $A_i$, and in terms of the measure $\mu$ are
given by
\begin{equation}
p_i = \frac{\mu(A_0 \cap L^{-1} (A_i))}{\mu (A_0)}. \label{probtrans}
\end{equation}

It is evident in the definition of $\Pi$ the dependence on the
particular sequence $(\epsilon_k)$ that we consider to define the map
$L$ and on the probability measure $\mu$.

It has been shown by Gaspard and Wang (\cite{GW}), that a way to
estimate the AIC of a string obtained from our dynamical system is the
theory of {\it recurrent events} applied to Markov chains
(\cite{Feller2},\cite{Feller}). From the theory of Markov chains, we
have that our stochastic matrix $\Pi$ is {\it irreducible}, and that
the state $A_0$ is {\it persistent}, in the sense that
\[
p[Z_m = 0 \hbox{ for some } m >n | Z_n = 0] =1 \quad \forall \ n,
\]
where $p$ is the probability measure on $(\Sigma_M,{\cal C},T_M)$ that
makes it a Markov chain with $(Z_n)$ as random variables from
$\Sigma_M$ to $\N$ (see Theorems \ref{yonescotulcea},
\ref{compatibili}).

If we consider as recurrent event ${\cal E}$ the passage from the
sub-interval $A_0$, we have that ${\cal E}$ is {\it certain} and that
the {\it mean recurrence time} $m_0$ is given by 
\begin{equation}
m_0 = \sum_{k=1}^{+\infty} \ k \ p_{k-1}. \label{mrt}
\end{equation}

If $m_0 = +\infty$ the state $A_0$ is called {\it null}, otherwise it
is {\it ergodic}. Thanks to the irreducibility of the stochastic
matrix $\Pi$, we have that all the states $A_i$ are of the same kind
of $A_0$, and then either $m_i = +\infty$ for all $i \in \N$ or $m_i
< +\infty$ for all $i \in \N$. 

For the recurrent event ${\cal E}$ two random variables can be
introduced: $X_k : \Sigma_M \to \N$ given by 1 plus the number of
trials between the $(k-1)$-th and $k$-th occurrence of ${\cal E}$;
$N_k : \Sigma_M \to \N$ given by the number of realizations of ${\cal
E}$ in $k$ trials. The random variables $X_k$ have all the same
probability distribution given by
\[
p[X_k = r] = p_{r-1},
\]
and their mean $E_p[X_k]=m_0$, where the subscript $p$ specifies the
measure we use to find the mean. For our problem it will be very
important also the form of the distribution function $F(x) =
\sum_{r=0}^{[x]} p_r$ of the $X_k$. Finally, we consider also the
probabilities $u_n$ that ${\cal E}$ occurs at the $n$-th trial. It
holds that
\begin{equation}
\lim_{n\to +\infty} u_n = \ \frac{1}{m_0}. \label{limitediu}
\end{equation}

Let's now explain the link between the AIC of a string generated by
the map $L$ and the theory of recurrent events (\cite{GW}). Given an
initial point $x\in I$ we obtain a string $\sigma \in \Sigma_M$ such
that $L^k(x) \in A_{\sigma_k}$ for all $k \in \N$. The string $\sigma$
is, for example, of the form
\begin{equation}
\sigma = (7654321054321002103210\dots). \label{stringa}
\end{equation}

One possible way to give an estimate for the AIC of the string is to
consider a compression of the string, and study the binary length of
the compressed string. One possible compression of the string $\sigma$
is given by
\begin{equation}
S = (75023\dots), \label{stringacompressa}
\end{equation}
that is the sequence of recurrent times for ${\cal E}$. So to the
finite string $\sigma^n$, obtained by the first $n$ symbols of
$\sigma$, we associate the string 
\[
S^{N_n}= (\sigma_{q_{_1}} \sigma_{q_{_2}} \dots \sigma_{q_{_{N_n}}})
\]
with $\sigma_{q_{_i} -1} =0$ for all $i$, where $N_n$ is the number of
realizations of ${\cal E}$. Then to have an idea of the behaviour of
the AIC of a string we have to estimate the behaviour of the random
variables $N_n$.

In \cite{Feller}, some possible behaviours for $E_p[N_n]$ have been
studied, for particular forms of the distribution function $F(x)$. In
particular

\begin{teorema}[Feller]
If the recurrence time of ${\cal E}$ has finite mean $m_0$ and
variance $V$, then
\[
E_p[N_n] \sim \ \frac{n}{m_0} + \frac{V-m_0+m_0^2}{2m_0^2}.
\]
If, instead $V=+\infty$, and the distribution function $F(x)$ satisfies
\[
F(x) \sim 1- A x^{- \alpha}
\]
with a constant $A$ and $0 < \alpha <2$, then:

i) If $1<\alpha <2$,
\[
E_p[N_n] \sim \ \frac{n}{m_0} + \frac{A}{(\alpha -1)(2-\alpha) m_0^2}
\ n^{(2-\alpha)};
\]

ii) If $0<\alpha <1$,
\[
E_p[N_n] \sim \ \frac{\sin \alpha \pi}{A \alpha \pi} n^{\alpha}.
\] \label{teofeller}
\end{teorema}

If the variance $V$ of the recurrence time is infinite but the
distribution function has a form different from that studied in
Theorem \ref{teofeller}, then we can show that

\begin{teorema}
If the state $A_0$ is ergodic then $E_p[N_n] \sim kn$, if instead
$A_0$ is a null state then $E_p[N_n]$ is an infinite of order
less than $n$.
\label{teoinf}
\end{teorema}

\noindent {\bf Proof.} The proof is based on a characterization of the
mean $E_p[N_n]$. In \cite{Feller}, it is shown that $E_p[N_n] = U_n -1$,
where
\[
U_n = \sum_{i=0}^n \ u_i.
\]

Then it follows that
\[
\lim_{n\to +\infty} \frac{E_p[N_n]}{n} = \lim_{n\to +\infty}
\frac{U_n}{n} = \lim_{n\to +\infty} u_n = \ \frac{1}{m_0}.
\]

Then, if $m_0<+\infty$, that is $A_0$ is ergodic, then $E_p[N_n]$ is
linear on $n$. Whereas if $m_0 = +\infty$, that is $A_0$ is a null
state, then $E_p[N_n]$ is an infinite of order less than $n$. \qed

\vskip 0.5cm The AIC calculated with the compression we have chosen
can be linked to the random variables $N_n$ by the following theorem

\begin{teorema}
For any string $\sigma \in \Sigma_M$ it holds
\begin{equation}
(N_n -1) + \log_2 (n-N_n +2) \leq I_{AIC}(\sigma^n) \leq N_n \log_2 \left(
\frac{n+N_n}{N_n} \right)
\label{eqstimainf}
\end{equation} 
\label{teostimainf}
\end{teorema}

\noindent {\bf Proof.} Given $n \in \N$, we have that $\Sigma_M = 
\cup C^n_r$, where the union is made on all the possible cylinders 
$C^n_r$, with $r \in \N^n$. Moreover we write 
\[
I_{AIC}(\sigma^n) = \sum_{i=1}^{N_n} \log_2 (\sigma_{q_{_i}} +2)
\]
for the AIC, where $\sigma_{q_{_i}} +2$ is used instead of
$\sigma_{q_{_i}}$, to have $\log_2 (\sigma_{q_{_i}} +2) \geq 1$ for
all $i$.

First of all, we consider only the cylinders $C^n_r$ with $r_n
=0$. This is done because we want to study the strings whose
compression changes when we increase our given $n$. Indeed, if the
compression wouldn't change, we wouldn't have any hint on the
behaviour of the AIC with respect to the length of the string.

We start with some special cases. Let's consider first the case $N_n
=n$. The only possible cylinder is then $C^n_r$ with
$r=(0,\dots,0)$. Then our compression doesn't change any string in
this cylinder, and $I_{AIC}(\sigma^n)= n $ for any string.

In the case $N_n =(n-1)$, the possible cylinders are given by $r \in
\N^n$ with only one symbol $1$ and all the others $0$. For all the
strings in these cylinders, the compression is $(n-1)$ symbols long,
and $I_{AIC}(\sigma^n) = (n-2)+\log_2 3$.

In the case $N_n =1$, the only possible cylinder is given by $r_i =
(n-i)$, and for the strings in this cylinder $I_{AIC}(\sigma^n) = \log_2
(n+1)$.

Let now be in general $N_n = n-h$, for some $h < n$. The
compression of such strings is then $N_n$-symbols
long. Moreover the compression is such that $\sum_{i=1}^{N_n}
\sigma_{q_{_i}} = h$. We now want to find the maximum and the minimum of
the function
\[
\sum_{i=1}^{N_n} \log_2 (\sigma_{q_{_i}} +2)
\]
with the condition $\sum_{i=1}^{N_n} \sigma_{q_{_i}} = k$. The maximum
is attained for equal $\sigma_{q_{_i}} \not= 0$, and the minimum for
all the $\sigma_{q_{_i}} =0$ but one which is equal to $h$.  Then the
maximum is given by $\sigma_{q_{_i}} = \frac{h}{n-h}$ for all $i$, and
the AIC for the strings in such a cylinder is given by
$\sum_{i=1}^{n-h} \log_2 \left( \frac{h}{n-h} +2 \right) = N_n \log_2
\left( \frac{n+N_n}{N_n} \right)$, and the minimum is given by
$(n-h-1) + \log_2 (h+2) = (N_n -1) + \log_2 (n-N_n +2)$.

Looking back at the special cases we studied before, we see that for
$N_n=n$ and $N_n=1$, the maximum and the minimum are the same, and
give exactly the value of the AIC we found building the
sequences. For the case $N_n = (n-1)$, the only possible value for the
AIC is equal to the minimum we found. This shows that
the maximum is not attained always, but there are cases in
which it is attained. For example let $n=8$, $N_n=4$, and consider the
string $\sigma = (10101010)$. Its compression is then given by
$S=(1111)$, and $I_{AIC}(\sigma^n) = 4\log_2 3 = N_n \log_2 \left(
\frac{n+N_n}{N_n} \right)$.

Finally we remark that we have tacitly assumed that our strings do not
begin with the symbol $0$. If it happened, the estimates wouldn't
change significantly.

We have thus proved that for a subset of $\Sigma_M$ of full measure
with respect to all the cylinders with $0$ as last symbol, the
AIC of a string can be estimated using the value of the random
variables $N_n$. \qed

\vskip 0.5cm Then our plan is the following: given a probability
measure $\mu$ on $(I,{\cal B},L)$, we find the stochastic matrix
$\Pi$, equation (\ref{stochmatrix}), the distribution function $F(x)$,
and, the mean $E_p[N_n]$. Then we can link the mean
$E_{\mu}[I_{AIC}(\sigma^n)]$, where $\mu$ is the probability measure on
$(I,{\cal B})$, with the mean $E_p[N_n]$, by Theorem
\ref{teostimainf}.

Another important aspect is the existence of an invariant measure for
the Markov chain associated to our dynamical system. Given a
stochastic matrix $\Pi$ of the form of equation (\ref{stochmatrix}) we
have 

\begin{teorema}
There is a measure $\bar p$ on the space $(\Sigma_M,{\cal
C},T_M)$, invariant for the stochastic matrix $\Pi$, defined by
\[
\bar p (k) = \sum_{n=0}^{+\infty} \ p_{k+n}.
\]
This measure is a probability measure if and only if the mean
recurrence time $m_0$ is finite.
\end{teorema}

This theorem allows us to induce on $(I,{\cal B},L)$ an $L$-invariant
measure $\bar \mu$, which is a probability measure if and only if
$A_0$ is ergodic.

\subsection{Restriction to the Lebesgue measure} \label{ssRttlm}
In the previous subsection we have shown that given a probability
measure $\mu$ on the space $(I,{\cal B},L)$, we have that the mean of
the AIC behaves differently according to the mean recurrence time of
the passage for the sub-interval $A_0$. These results clearly depend
on the choice of the measure $\mu$ and of the sequence $(\epsilon_k)$
used in the definition of the map $L$. In this subsection we want to
study these two problems from a practical point of view.

When we apply the notion of AIC to a string obtained from a dynamical
system, the choice of this string depends on the choice of the initial
point $x$ which we use to generate the orbit of the dynamical
system. This choice can be made randomly, and the most natural way to
introduce a probability distribution on the choice of the initial
point is by using the Lebesgue measure $l$ on the space. Hence we
apply the results of Subsection \ref{ssGr} to the system $(I,{\cal
B},L,l)$. Using the Lebesgue measure $l$, thanks to the piecewise
linearity of the map $L$, the probabilities of transition $p_i$ given
in equation (\ref{probtrans}) assume a particular simple form. Indeed
we find that $p_i = l(A_i)$ for all $i \in \N$.

With respect to the Lebesgue measure $l$, it is also possible to prove
that the compression we have introduced for strings given by any map
$L$ is the best possible. We have

\begin{teorema}
Given any piecewise linear map $L$ of the form (\ref{mannmaplin}), the
best compression for the strings generated by the dynamical system
$(I,{\cal B},L,l)$, where $l$ is the Lebesgue measure, is the
compression given in equations (\ref{stringa}) and
(\ref{stringacompressa}), for $l$-almost any initial condition.
\label{bestcompressiontheorem}
\end{teorema}

\noindent {\bf Proof.} The compression we introduced gives a bijective
relation between our space $(\Sigma_M,{\cal C},T_M,p)$ and the space
of all possible infinite sequences built on countable symbols, without
any restriction given by a transition matrix. We denote this space as
$\Sigma$. We introduce on $\Sigma$ the $\sigma$-algebra of the
cylinders $C^{'n}_r$, given as in equation (\ref{cilindri}), and a
probability measure $p'$ inherited in some way from $p$. We define
$p'$ by
\begin{equation}
p'(C^{'n}_r) = p(C^N_R), \label{confmis}
\end{equation}
where the cylinder $C^N_R$ is built in such a way that compression of
strings that belong to it gives strings belonging to the cylinder
$C^{'n}_r$. At this point we ask if it is possible to compress any
more strings belonging to the space $(\Sigma,{\cal C'},p')$. But, if
on this space we consider the usual shift map, we find a
$p'$-invariant map with positive K-S entropy $h^{KS}$. This is given
by a direct computation, using the piecewise linearity of the map $L$
that gives a simple form for the measure $p'$. Then it is well known
that the AIC for this dynamical system behaves like $h^{KS} n$, for
$p'$-almost any string $\sigma$, then for $l$-almost any initial
condition $x \in I$. This clearly implies that for $l$-almost any
initial condition $x \in I$, ours is the best possible
compression. \qed

\vskip 0.5cm
\begin{guess}
According, for example, to Brudno's approach (\cite{Brudno}) to obtain
a definition of AIC for finite orbits of a dynamical system, we should
evaluate the supremum of $I_{AIC}(\sigma^n)$ varying the open covers
of the interval $I=[0,1]$. Theorem \ref{eqinfdinsimb} can be
established as well if we consider an open cover of the form
$A_i=(\epsilon_i, \epsilon_{i-1} +\eta_i)$ for $0<\eta_i <<1$. Theorem
\ref{bestcompressiontheorem} suggests that the AIC of any sequence
generated by the system $(I,{\cal B},L)$ with a non-trivial open cover
has to be related to the random variables $N_n$ as in the case we are
considering. Then we can prove that the AIC of the particular strings
we are considering is the AIC of the dynamical system $(I,L)$. This
can also be proved in a more general contest (\cite{Stefano}).
\end{guess}

The second point is the choice of the sequence $(\epsilon_k)$. We know
that for any sequence, we can find an isomorphism of the map
$L$ with the Manneville map $f$ of equation (\ref{mannmap}). At this
point we present some cases for the choice of $(\epsilon_k)$, and then
study a particular choice. 

\begin{esempio}
Let the sequence be given by $\epsilon_k = \frac{1}{k^{\alpha}}$ with
$\alpha >0$. This sequence has all the properties we need for the
definition of the map $L$. Then
\[
p_i = l(A_i) = \ \frac{1}{(i-1)^\alpha} - \frac{1}{i^\alpha} \sim \
\frac{1}{i^{\alpha +1}}.
\]

The mean recurrence time $m_0$ is given by
\[
m_0 = \sum_{r=1}^{+\infty} \ r \ p_{r-1}
\]
and it is finite if and only if $\alpha >1$. Then we can find an
invariant measure $\bar \mu$ for the system $(I,{\cal B},L)$, such
that $\bar \mu (A_i) \sim \frac{1}{i^\alpha}$. This measure $\bar \mu$
is a probability measure if and only if $\alpha >1$.

The variance $V$ of the recurrence time is given by
\[
V = \sum_{r=1}^{+\infty} \ r^2 \ p_{r-1}
\]
and then $V<+\infty$ if and only if $\alpha > 2$. If we compute the
distribution function $F(x)$, we have that
\[
F(x) \sim 1- A x^{-\alpha},
\]
then we can apply Theorems \ref{teofeller} and \ref{teostimainf}, and
obtain that:
\begin{itemize}
\item[i)] if $\alpha > 1$ then $E_l[I_{AIC}(\sigma^n)] \sim E_p[N_n]
\sim n$;

\item[ii)] if $\alpha <1$ then $E_p[N_n] \sim n^\alpha$ and $n^\alpha
\leq E_l[I_{AIC}(\sigma^n)] \leq n^\alpha \log_2 n$ (see Theorem
\ref{teostimainf}).
\end{itemize}
\label{esempio1}
\end{esempio}

\begin{esempio}
Let now the sequence $(\epsilon_k)$ be given by $\epsilon_k =
\frac{1}{a^k}$, where $a \in \N$ and $a >1$. In this case it is easy
to verify that $p_i \sim \frac{1}{a^i}$. Then the mean recurrence time
$m_0$ and the variance $V$ of the recurrence time are always finite,
and we can deduce that $E_l[I_{AIC}(\sigma^n)] \sim E_p[N_n] \sim n$. The
invariant probability measure $\bar \mu$ exists and is given by $\bar
\mu (A_i) \sim \frac{1}{a^i}$.
\label{esempio2}
\end{esempio}

\begin{esempio}
Finally let's consider the case in which the sequence $(\epsilon_k)$
is given by $\epsilon_k = \frac{1}{k^\alpha (\log k)^\beta}$ with
either $\alpha =1$ and $\beta >1$ or $\alpha >1$. If we compute the
variance $V$ of the recurrence time, we obtain that it is finite if
and only if either $\alpha >3$ or $\alpha=3$ and $\beta>1$. In this
case we obtain from Theorem \ref{teofeller} that
$E_l[I_{AIC}(\sigma^n)] \sim n$. But if $V=+\infty$, we cannot find an
explicit form for the distribution function $F(x)$ similar to that of
Theorem \ref{teofeller}, but we can use Theorem \ref{teoinf}. Indeed
we have that the state $A_0$ is ergodic if and only if either $3
>\alpha >2$ or $\alpha =2$ and $\beta >1$, in which case
$E_l[I_{AIC}(\sigma^n)] \sim n$. For other values of $\alpha$ and
$\beta$, we know that $E_p[N_n]$ is an infinite of order less than
$n$, and for $E_l[I_{AIC}(\sigma^n)]$ we can apply Theorem
\ref{teostimainf} to obtain an estimate for its order of infinite.
\label{esempio3}
\end{esempio}

\begin{esempio}
In a particular case, we can say something more about the order of
infinite of $E_p[N_n]$ using the theory of recurrent events
(\cite{Feller}) and of power series (\cite{Titchmarsh}). Indeed,
choosing the sequence $\epsilon_k \sim \frac{1}{\log k}$, we have that
the distribution function $F(x) \sim \left( 1 - \frac{1}{\log x}
\right)$. Then from the characterization of $E_p[N_n]$ in terms of the
generating function of the random variables $X_k$, we obtain that
asymptotically $E_p[N_n] < n^\alpha$ for all $\alpha >0$. Hence, from
Theorem \ref{teostimainf}, $E_l[I_{AIC}(\sigma^n)]$ is an infinite of
order smaller than any power law.
\label{esempio4}
\end{esempio}

We have thus found that changing the sequence $(\epsilon_k)$ it is
possible to obtain all a set of behaviours for the mean of the AIC
with respect to the Lebesgue measure. It is clear that the different
behaviours are induced by the rate with which the derivative of the
map $L$ increases, rate that depends on the order with which
$(\epsilon_k)$ tends to 0. It is then evident that this must be also
the criterion to distinguish between different behaviours of the
information function of the Manneville map $f$ for different values of
the parameter $z$ in equation (\ref{mannmap}). We have then to find a
way to associate to a particular value of $z$ a given sequence
$(\epsilon_k)$. Since we have to maintain a given rate of increasing
of the derivative, given a value for $z$, we look for the sequence
$(\epsilon_k)$ such that $\epsilon_k \sim x_k$, where $x_k$ is the
sequence of preimages of the point $x_0$, as defined in Section
\ref{sTlm}.

We have that $x_k \sim \frac{1}{k^\alpha}$ with $\alpha =
\frac{1}{z-1}$. We are then in the first case we studied, and $\alpha
>0$ being $z >1$. We can apply all the results we found, in particular:
\begin{itemize}
\item if $z <2$ then $E_l[I_{AIC}(\sigma^n)] \sim n$ and there is an
$L$-invariant probability measure $\bar \mu$ such that $\bar \mu (A_i)
\sim \frac{1}{i^\alpha}$;

\item if $z>2$ then $n^\alpha \leq E_l[I_{AIC}(\sigma^n)] \leq n^\alpha
\log_2 n$ and the invariant measure $\bar \mu$ is not a probability
measure.
\end{itemize}

We have then found as a particular case the same results of
\cite{GW}. But, we would like to have a behaviour of the AIC valid for
almost any orbit with respect to the Lebesgue measure $l$. We have

\begin{teorema}
For almost any point $x \in I$ with respect to the Lebesgue measure
$l$ and for all $\delta>0$, we have that
$E_{\mu_\delta}[I_{AIC}(\sigma^n)]$ is asymptotically equivalent to
$E_l[I_{AIC}(\sigma^n)]$, where with $\mu_\delta$ we denote the
measure given by the Lebesgue measure $l$ concentrated on
$U_\delta=(x-\delta,x+\delta)$.
\label{teoqo}
\end{teorema}

\noindent {\bf Proof.} The proof is based on a simple application of
the method we used before. Indeed, given $x \in I$, let's consider
$U_\delta=(x-\delta,x+\delta)$ for some $\delta$. If we now want to
estimate the value of $E_\delta[N_n]$, that is the mean made with
respect to the measure on $\Sigma_M$ induced by $\mu_\delta$, we
notice that, from the properties of the map $L$, it is clear that
there exists a $R(\delta)\in \N$ such that $L^{R(\delta)} (U_\delta) =
[0,1]$. We can deduce that $E_\delta[N_n]$ is the same as $E_p[N_n]$,
where $p$ is the measure on $\Sigma_M$ induced by $l$, for
$n>R(\delta)$. Then it follows that $E_{\mu_\delta}[I_{AIC}(\sigma^n)]
\sim E_l[I_{AIC}(\sigma^n)]$ for $n>R(\delta)$.  \qed

\begin{cor}
The AIC of the Manneville map $f$ is given by
\begin{itemize}
\item $n^\alpha \leq E_{\mu_\delta}[I_{AIC}(\sigma^n)] \leq n^{\alpha}
\log_2 n$ with $\alpha = \frac{1}{z-1}$, for $z>2$;

\item $E_{\mu_\delta}[I_{AIC}(\sigma^n)] \sim n$ for $z<2$
\end{itemize}
for almost any point $x\in I$ with respect to the Lebesgue measure
$l$. \label{cor2}
\end{cor}

We remark that there are points $x\in I$ for which $I_{AIC}(x^n) \sim
n$ also for the Manneville map $f$ with $z>2$. Indeed, in Subsection
\ref{ssTame}, we proved the existence of $f$-invariant measures $\mu$
on $(I,{\cal B})$, for which the K-S entropy $h^{KS}_\mu =\log
2$. This is not a contradiction since such measures have support on a
set of zero Lebesgue measure.

\section{Conclusions}

In this paper we have proved that the Manneville map with $z>2$
(\ref{mannmap}) exhibits, from the AIC point of view, a behaviour
which is intermediate between the so-called {\it full chaos} (positive
K-S entropy) and periodicity. Then we obtain that the complexity (see
Section \ref{sI}, equations (\ref{complessita})) of the Manneville map
with $z>2$ is null for almost any initial condition with respect to
the Lebesgue measure $l$ and then the algorithmic entropy $h_l =0$
(see Section \ref{sI}, equations (\ref{entropia})). In particular we
have found a family of piecewise linear maps $L$ that, for fixed
sequences $\epsilon_k$ (see equation (\ref{mannmaplin})), can be used
as a model for the Manneville map. Moreover this family of maps
presents a rich set of possible algorithmic behaviours (depending on
the choice of the map $L$ of the family). It is evident that changing
the sequence $\epsilon_k$, the algorithmic behaviour varies from full
chaos to {\it mild chaos}, which is characterized by
$I_{AIC}(\sigma^n)$ of order smaller than any power law. This
behaviour can be achieved from the Manneville map at the limit $z \to
+\infty$ (see Theorem \ref{teostimainf}), and Example \ref{esempio4}
seems to suggest a way to find a particular sequence generating mild
chaos. We remark that it would be impossible using the notions of
complexity and algorithmic entropy to distinguish many of the maps of
the family $L$, since we would obtain $K(x)=h_l =0$ for almost every
$x \in (I,{\cal B},l)$. Then the AIC is for those maps a powerful tool
for the classification and it can be actually estimated.

Finally we observe that it has been proved that the AIC of a string is
not a computable function and, in particular, it cannot be computed by
any algorithm (\cite{Chaitin}). Then, apart from particular dynamical
systems for which we can estimate the AIC, the classification of
dynamical systems using the AIC cannot be obtained
explicitly. Nevertheless we believe that it is fundamental to obtain
an explicit estimation of the AIC for as many dynamical systems as
possible. However the AIC can be approximated by different notions of
information content of a string. In particular we can define
\begin{equation}
I_A(\sigma^n) = |A(\sigma^n)|, 
\end{equation}
the {\it information function}, where $|A(\sigma^n)|$ is the binary
length of the output obtained from a string $\sigma^n$ by means of a
compression algorithm $A$ (see \cite{Licata}). A particular
compression algorithm, called {\it CASToRe}, has been built to analyze
dynamical systems which present rich dynamics, but having zero K-S
entropy for all invariant measures that are physically relevant
(\cite{Argenti}). The algorithm {\it CASToRe} has been tested on the
Manneville map, giving as result $I_{CASToRe}(\sigma^n)\sim n^\alpha$
for $\alpha <1$, confirming our results (\cite{Argenti},
\cite{Licata}), and on the logistic map at the chaos threshold, giving
the presence of mild chaos (\cite{Bonanno}). At the moment it is not
clear whether the algorithm {\it CASToRe} can approximate the AIC for
any dynamical system, but from this paper on the Manneville map and
from many experimental results on the logistic map, it seems that at
least in these two cases there is evidence of accordance between the
theoretical predictions and the experiments with {\it CASToRe}.

\end{document}